\documentclass[12pt]{article}
\usepackage{amssymb}

\usepackage{amsthm}
\usepackage{amsmath}

\newtheorem{defi}{\bf Definition}[section]
\newtheorem{theo}[defi]{\bf Theorem}

\newtheorem{coro}[defi]{\bf Corollary}

\newtheorem{lem}[defi]{\bf Lemma}

\title{Lie structure in semiprime superalgebras with superinvolution}

\author{Jes\'us Laliena  \footnote {The first author has been
supported by the Spanish Ministerio de Educaci\'on y Ciencia (MTM 2004-08115-CO4-02) and both by the Comunidad Aut\'onoma de La Rioja (ANGI 2005/05).} and Sara Sacrist\'an
\\{\small Departamento de Matem\'aticas y Computaci\'on}\\
{\small  Universidad de La Rioja}\\
{\small  26004, Logro\~no. Spain}\\
{\small jesus.laliena@dmc.unirioja.es \quad ssacrist@ya.com}}
\date{\quad}

\begin{document}
\maketitle\vspace{-1.5cm}

\begin{abstract}

In this paper we investigate the Lie structure of the Lie superalgebra $K$ of skew elements of a
semiprime associative superalgebra $A$ with superinvolution. We show that if $U$ is a Lie ideal
of $K$, then either there exists an ideal $J$ of $A$ such that the Lie ideal   $[J\cap K,K]$ is nonzero and contained in $U$, or  $A$ is a subdirect sum of $A'$, $A''$, where the image of $U$ in $A'$
is central, and $A''$ is a subdirect product of orders in simple superalgebras, each at most
16-dimensional over its center.

\end{abstract}

{\parindent= 4em \small  \sl Keywords: associative superalgebras, semiprime superalgebras, superin-}

{\parindent=10em \small  \sl  volutions, skewsymmetric elements, Lie structure.}

\section{Introduction.}

\bigskip

The study of the relationship between the structure of an associative algebra $A$ and that of the Lie
algebra $A^-$ was  started by I. N. Herstein (see \cite {H1}, \cite{H2}) and W. E. Baxter (see \cite {B}). Afterwards, several authors have made different contributions and generalizations to the subject (see for instance  \cite{Er}, \cite {La}, \cite {Ma} ). 

Regarding superalgebras, this line of research was motivated by the classification of the finite dimensional simple Lie superalgebras given   by V. Kac (\cite {Ka}), particularly the types given from simple associative superalgebras and from simple associative superalgebras with superinvolution.  In \cite {Go-S}, thinking in simple associative superalgebras with superinvolution, C. G\'omez-Ambrosi and I. Shestakov investigated the Lie structure of the set of skew elements, $K$, of a simple associative superalgebra, $A$, with superinvolution over a field of characteristic not 2. These results were later extended to prime associative
superalgebras with superinvolution (\cite {Go-L-S}). It was specifically proved that the Lie ideals
of $K$ and $[K,K]$ are of the kind $[J\cap K, K]$ for a nonzero ideal $J$ of $A$, if $A$ is
nontrivial, that is  with a nonzero odd part, and if $A$ is not a central order in a Clifford
superalgebra with at most 4 generators.

This paper is devoted to the description of the Lie ideals of $K$, the set of skew elements of  a
semiprime associative superalgebra, $A$, with superinvolution * over a commutative unital ring
$\phi$ of scalars with ${1\over 2 } \in \phi$.

We notice that the Lie structure of prime superalgebras and simple superalgebras without superinvolution was studied by F. Montaner (see \cite {M}) and S. Montgomery (see \cite {Mo}).

For a complete introduction to the basic definitions and examples of superalgebras,
superinvolutions and prime and semiprime superalgebras, we refer the reader to \cite {Go-S} and \cite
{M}.

Throughout the paper, unless otherwise stated,  $A$ will denote a nontrivial semi-prime associative
superalgebra with superinvolution * over a commutative unital ring $\phi$ of scalars with ${1\over
2}\in \phi$. By a nontrivial superalgebra we understand a superalgebra with nonzero odd part. $Z$ will
denote the even part of the center of $A$, $H$ the Jordan superalgebra of symmetric elements of $A$,
and $K$ the Lie superalgebra of skew elements of $A$. If $P$ is a subset of $A$, we will denote by
$P_H=P\cap H$ and $P_K= P\cap K$. The following containments are straightforward to check, and they
will be used throughout without explicit mention:
$[K,K]\subseteq K, \quad [K,H]\subseteq H, \quad [H,H]\subseteq K, \quad H\circ H \subseteq H, \quad
H\circ K\subseteq K$ and
$K\circ K\subseteq H$.

We recall that a superinvolution * is said to be of the first kind if $Z_H=Z$, and 
of the second kind if $Z_H\not= Z$.

If $Z\not= 0$, one can consider the localization $Z^{-1}A=\{z^{-1}a : 0\not= z \in Z, a\in A\}$. If
$A$ is prime, then  $Z^{-1}A$ is a central prime associative superalgebra over the field
$Z^{-1}Z$.  We call this superalgebra  the central closure of $A$. We also say that $A$ is a
central order in $Z^{-1}A$. While this terminology is not the standard one, for which the definition involves
 the extended centroid,  if $Z\not=0$ both notions coincide (for more specifications see 1.6 in \cite
{M}). 

Let $A$ be a prime superalgebra, and let  $V=Z_H-\{0\}$ be the subset of regular symmetric elements. Note that if $Z\not= 0$, $Z_H\not=0$. Also  $Z^{-1}A=V^{-1}A$, since for all  $0\not=z \in Z, a\in A$ we have $z^{-1}a=(zz^*)^{-1}(z^*a)$. It will be more convenient for us, in order to extend the superinvolution in a natural way, to work with $V$ rather than with $Z$.  We may consider $V^{-1}A$ as a superalgebra over the field $V^{-1}Z_H$. Then the superinvolution on $A$ is extended to a superinvolution  of the same kind on $V^{-1}A$ over  $V^{-1}Z_H$
via $(v^{-1}a)^*=v^{-1}a^*$. It is then easy to check that $H(V^{-1}A,*)=V^{-1}H$ and
$K(V^{-1}A,*)=V^{-1}K$. Moreover, $Z(V^{-1}A)_0=V^{-1}Z$ and $V^{-1}Z\cap V^{-1}H=V^{-1}Z_H$. We will say that the superalgebra $V^{-1}A$ over the field $V^{-1}Z_H$ is the *-central closure of $A$.

 We notice that in every semiprime superalgebra $A$,  the intersection of all the prime ideals $P$
of $A$ is zero. Consequently $A$ is a subdirect product of its prime images. If each prime image of
$A$ is a central order in a simple superalgebra at most $n^2$ dimensional over its center, we  say
that $A$ verifies $S(n)$.

If $M$ is a subsupermodule of $A$, we denote by $\bar M$ the subalgebra of $A$ generated by $M$. We will say that $M$ is dense if $\bar M$ contains a nonzero ideal of $A$.

In this paper, we prove that if $K$ is the Lie superalgebra of skew elements of a semiprime
associative superalgebra with superinvolution, $A$, and $U$ is a Lie ideal of $K$, then one of the following alternatives must hold: either $U$ must contain a nonzero Lie ideal $[J\cap K,K]$, for $J$ an ideal of $A$, or  $A$ is a subdirect sum of $A^\prime$, $A^{\prime \prime}$, where the image of $U$ in $A^{\prime \prime}$ is central and $A^\prime $ satisfies $S(4)$.

The following results are instrumental for the paper:

\begin{lem} (\cite {H2}, lemma 1.1.9)
If $A$ is a semiprime algebra and $[a,[a,A]]=0$, then $a\in Z(A)$.
\end{lem}

\begin{lem} (\cite {M}, lemmata 1.2, 1.3)
If $A=A_0 \oplus A_1 $ is a prime superalgebra, then $A$ and $A_0$ are semiprime and
either $A$ is prime or $A_0$ is prime (as algebras).
\end{lem}

\begin{lem} (\cite {M}, lemma 1.8)
Let $A=A_0 \oplus A_1$ be a prime superalgebra. Then

\begin{enumerate}
\item[{\rm (i)}] If $x_1 \in A_1 $ centralizes a nonzero ideal $I$ of $A_0$, then
$x_1 \in Z(A)$.

\item[{\rm (ii)}] If $x_1 ^2$ belongs to the center of a nonzero ideal $I$ of $A_0$, then
$x_1 ^2 \in Z(A)$.
\end{enumerate}
\end{lem}

\bigskip

\begin{lem} (\cite {Go-L-S}, Corollary 2)
Let $A$ be a semiprime superalgebra and $L$ a Lie ideal of $A$. Then either $[L,L]=0$, or $L$ is dense
in $A$.
\end{lem}

\begin{lem}(\cite {Go-L-S}, Theorem 2.1)
Let $A$ be a prime nontrivial associative superalgebra. If $L$ is a Lie ideal of $A$,  then either
$L\subseteq Z$ or $L$ is dense in $A$, except if $A$ is a central order in a 4-dimensional Clifford
superalgebra.
\end{lem}

We remark that the bracket product in lemma 1.1 is the usual one: $[a,b]=ab-ba$, but the bracket
product in lemmata 1.3,1.4,1.5 is the superbracket $[x_i,y_j]_s=x_iy_j-(-1)^{ij}y_jx_i$ for $x_i\in
A_i, y_j\in A_j$ homogenous elements. In fact, the superbracket product coincides with the usual
bracket if one of the arguments belongs to the even part of $A$. In the following, to simplify the
notation, we will denote both in the usual way $[\ , \ ]$ but we will understand that it is the
superbracket if we are in a superalgebra.

\bigskip

\section {Lie structure of K.}

\bigskip

Let $A$ be an associative superalgebra and $M,S$ be subgroups of $A$. Define $(M:S)=\{a\in A :
aS\subseteq M\}$.

 Let $U$ be a Lie ideal of $K$. We recall (see lemma 4.1 in \cite {Go-S}) that $K^2$ is a Lie ideal of
$A$.

\begin{lem}
If $A$ is semiprime, then either $U$ is dense in $A$ or $[u\circ v,w]=0$
for every $u,v,w\in U$.
\end{lem}

\begin{proof}[Proof:]
We have 
 $$[u\circ v, k]= u\circ [v,k] + (-1)^{{\bar k}{\bar v}} [u,k]\circ v\in \bar U$$
{\noindent for every
$u,v \in U$ and $k\in K$. And also for any $u,v\in U$ and $h\in H$ we get}
$$[u\circ v,h]=[u,v\circ h]+(-1)^{\bar u \bar v}[v,u\circ h]\in  U,$$
{\noindent because $K\circ H\subseteq K$. Since $A= H\oplus K$ it
follows that  $[u\circ v, A]\subseteq \bar U$ for any
$u,v\in U$. But for any $a\in A$}
$$ [u\circ v, wa]=[u\circ v, w]a + (-1)^{({\bar u}+{\bar v}){\bar w}} w[u\circ v,a]$$
{\noindent and so $[u\circ v,w] A\subseteq \bar U$ for every $u,v,w \in U$, that is, $[u\circ v,w]\in
 (\bar U:A)$. We notice that from the above equations we can also deduce that $[u\circ v,w]a\in
\bar {K^2}$ and so $[u\circ v,w]\in(\bar {K^2}:A)$.}

We claim that $A(\bar {K^2}:A)\subseteq (\bar {K^2}:A)$. Indeed, for any $x\in (\bar {K^2}:A), a,b\in
A$ 
$$ axb=(-1)^{(\bar x + \bar b)\bar a} (xb)a+[a,xb]\in \bar {K^2},$$
{\noindent because $[\bar {K^2},A]\subseteq \bar {K^2}$ (for any  $t,s\in K^2$, 
$[ts,a]=t[s,a]+(-1)^{\bar s \bar a}[t,a]s\in \bar K^2$). Hence
$A(\bar {K^2}:A)A\subseteq (\bar {K^2}:A)A\subseteq \bar {K^2}$. But $K(\bar U:A)\subseteq (\bar U:A)$
because for any $x\in (\bar U:A), k\in K, a\in A$}
$$(kx)a=[k,xa]+(-1)^{(\bar x +\bar a)\bar k}(xa)k\in \bar U,$$
{\noindent because $[K,\bar U]\subseteq \bar U$, and so $\bar {K^2}(\bar U:A)\subseteq (\bar U:A)$.
Therefore, we finally get} 
$$A(\bar {K^2}:A)A(\bar U:A)A\subseteq \bar {K^2}(\bar U:A)A\subseteq \bar U,$$
{\noindent and since $[u\circ v,w]\in (\bar U:A)$ and also $[u\circ v,w]\in (\bar {K^2}:A)$ it follows
that $A[u\circ v,w]A[u\circ v,w]A \subseteq \bar U$. Thus, since $A$ is semiprime, either $[u\circ
v,w]=0$ for any
$u,v,w\in U$ or $U$ is dense in $A$.}
\end{proof}

\bigskip

We note that the ideal contained in $\bar U$ in the above Lemma, $J=A[u\circ v,w]A[u\circ v,w]A$, is
also a
$*$-ideal, that is, $J^*\subseteq J$.

\begin{lem}
Let $A$ be semiprime, and let $U$ be a Lie ideal of $K$ such that $[U\circ U,U]=0$.  Then

\begin{enumerate}

\item[{\rm (i)}] $u\circ v\in Z$ for every $u,v \in U_0$.

\item[{\rm (ii)}] $u\circ v=0$ for every $u,v\in U_1$.
\end{enumerate}

\end {lem}

\begin{proof}[Proof:]
Assertion (i) is proved as in Theorem 5.3 of \cite {Go-S}, and  (ii)  as
in Theorem 3.2 of \cite {Go-L-S}.
\end{proof}

\bigskip

Next we deal with the second case of lemma 2.1, that is, when $[u\circ v,w]=0$ for any $u,v,w\in U$ (and
therefore when $u\circ v\in Z$ for every $u,v\in U_0,$ and $u\circ v=0$ for every $u,v\in U_1$),
and we will study the prime images of
$A$.

\bigskip

Let $P$ be a prime ideal of $A$. We will suppose first that $P^*\not= P$. In this case $(P^*+P)/P$ is
a nonzero proper ideal of $A/P$ and we claim that $(P^*+P)/P\subseteq (K+P)/P$. Indeed, if
$y\in P^*$ then $y+P=(y-y^*)+ y^*+P\in (K+P)/P$. Also if $U$ is a Lie ideal of
$K$ we have that $(U+P)/P$ is an abelian subgroup of $A/P$ and satisfies
$$[(U+P)/P,(P^*+P)/P]\subseteq ([U,K]+P)/P\subseteq (U+P)/P.$$

Therefore $(U+P)/P$ is a Lie ideal of $(P^*+P)/P$, and $(P^*+P)/P$ is an ideal in $A/P$, a prime
superalgebra. Of course if $u\circ v\in Z$ for every $u,v \in U_0$ and $u \circ v =0$ for any $u,v\in
U_1$, then the same property is satisfied in $A/P$, that is, $(u+P) \circ (v+P) \in Z_0 (A/P)$ for
every $u+P, v+P \in (U_0+P)/P$, and $(u+P)\circ (v+P)=0$ for any $u+P, v+P\in  (U_1+P)/P$. Let us analyze this situation.  We notice that the assumption that $A/P$ has a superinvolution is not required. We state first a useful lemma.

\begin{lem}
Let $A$ be a prime superalgebra, $I$ a nonzero ideal of $A$ and $U$ a subset of $A$ such that
$[U,I]=0$. Then $U\subseteq Z$.
\end{lem}

\begin{proof}[Proof:]
For any $u_k\in U_k, a_i\in A_i, y_j\in I_j$, applying $[U,I]=0$ we get
$$u_k(a_iy_j)= (-1)^{(i+j)k}
(a_iy_j)u_k=(-1)^{ik} a_i(u_ky_j).$$

{\noindent Since $A$ is prime it follows that $u_ka_i=(-1)^{ik}a_iu_k$. On the other hand, given
$u_1 \in U_1$ we have $[u_1,I_0]=0$, and applying lemma 1.3 (i), $u_1\in
Z_1 (A)$. Hence for every $u_1 \in U_1, a_1 \in A_1$ we have
$u_1 a_1 = a_1 u_1 = - a_1 u_1$, that is, $a_1 u_1 =0$, and, 
because $u_1 \in Z(A)$ and the primeness of $A$, $U_1 =0$ and $U\subseteq Z$.}

\end {proof}

\begin{theo}
Let $A$ be a prime superalgebra, and let $I$ be a nonzero proper ideal of $A$. Suppose that  $U$ is an abelian subgroup of 
$A$ such that $[U,I] \subseteq U$,  $u\circ v \in Z$ for every $u,v \in U_0$, and $u\circ v=0$ for
every $u,v\in U_1$. Then either $A$ is commutative, or $A$ is a central order in a 4-dimensional
simple superalgebra, or
$U\subseteq Z$.
\end{theo}

\begin{proof}[Proof:]

Let $T=\{x\in A: [x,A]\subseteq [U,I]\}$. Since
 $$[[[U,I],[U,I]],A]\subseteq [[U,I],[[U,I],A]]\subseteq [[U,I],I]\subseteq [U,I],$$

{\noindent we have $[[U,I],[U,I]]\subseteq T$. We notice that $T$ is subring because for any $t,s\in
T$,} 
$$[ts,a]=[t,sa]+(-1)^{\bar t \bar s + \bar a \bar t}
[s, at]\in [U,I].$$
{\noindent Let $T^\prime$ be the subring generated by $[[U,I],[U,I]]$. Since}
$$[[[U,I],[U,I]],I]  \subseteq [[U,I],[[U,I],I]] subseteq [[U,I],[U,I]]$$ 
{\noindent  it follows that
$[T^\prime,I]\subseteq T^\prime$. We consider now two cases: a) $[T^\prime,I]=0$, and b)
$[T^\prime,I]\not=0$.}

\bigskip

a) If  $[T^\prime,I]=0$, then $[[[U,I],[U,I]],I]=0$. By lemma 2.3 we get
$[[U,I],[U,I]]\subseteq Z$, and so $[[U,I],[U,I]]_1=0$.

We claim that in this situation either  $U\subseteq Z$, or $A$ is commutative, or $A$ is a central order in a
4-dimensional simple superalgebra. We present the proof of this in 6 steps.

\medskip

1. $[U,I]_0\subseteq Z$. By  hypothesis
$u\circ v\in Z$ for any
$u,v\in U_0$, so since
$[U,I]\subseteq U$ it follows that  $uv\in Z$ for any $u,v \in [U,I]_0$. Hence, for any $u,v \in
[U,I]_0$, we have 
$$[u,v][u,v]=[u,v[u,v]]-v[u,[u,v]]=[u,[vu,v]]-[u,[v,v]u]=0$$

{\noindent because $[u,v], vu\in Z$. Therefore, from the primeness of $A$, $[u,v]=0$ for any $u,v\in
[U,I]_0$. So since $[[U,I],[U,I]]_1=0$, $[u,[u,I]]=0$ for any $u\in [U,I]_0$, and therefore, by lemma
1.2 and theorem 1 in
\cite {H3},
$u\in Z(I)$, that is $[U,I]_0\subseteq Z$ because $A$ is prime.

\medskip

2. $[U_0, I_0]=0$. By step 1  we have
$[u_0,[u_0,I_0]]=0$ for any $u_0\in U_0$, and again by theorem 1 in \cite {H3} and lemma 1.2, we
obtain that
$[U_0,I_0]=0$. 

\medskip

3. $U_1 U_1\subseteq Z$.  Let $u_1\in U_1, y_1\in I_1$,  since $[U_1,I_1]\subseteq
[U,I]_0\subseteq Z$ we get
$$[u_1^2,y_1]=u_1[u_1,y_1]-[u_1,y_1]u_1=u_1[u_1,y_1]-u_1[u_1,y_1]=0.$$
Therefore, since $u_1\circ v_1=0$ for any $u_1, v_1\in U_1$, $0=[(u_1+v_1)^2,y_1]=[u_1v_1+v_1u_1,y_1]=
2[u_1v_1,y_1]$ for any $y_1\in I_1$. And, since $[u_1,v_1]=2u_1v_1\in U_0$ because $u_1\circ v_1=0$
for any $u_1,v_1\in U_1$, we have $[u_1v_1,I_0]=0$ for any $u_1,v_1\in U_1 $ by step 2. So $[u_1v_1,I]=0$ for
any $u_1,v_1\in U_1$, and then $u_1v_1\in Z$, because of lemma 2.3.}

\medskip

4. $I_1(U_1)^3\subseteq Z$. From the steps 1 and 3 for any $u_1,v_1,w_1\in U_1, y_1\in I_1$ we get
$[u_1,y_1]v_1w_1\in Z$, but

\begin{eqnarray*}
[u_1,y_1]v_1w_1&=&[u_1,y_1v_1]w_1+y_1[u_1,v_1]w_1\\&=&[u_1w_1,y_1v_1]-u_1[w_1,y_1v_1]+
y_1[u_1,v_1w_1]+y_1v_1[u_1,w_1]
\end{eqnarray*}

{\noindent and since $[u_1w_1,y_1v_1]=0$, $y_1[u_1,v_1w_1]=0$, by step 3 and
$u_1[w_1,y_1v_1]\in U_1[U_1,I_0]\subseteq
U_1U_1\subseteq Z$, we obtain that $y_1v_1[u_1,w_1]\in Z$, that is $I_1(U_1)^3\subseteq Z$ because
$u_1
\circ w_1=0$.}

\medskip

5. Either $U_1U_1=0$ or $A$ is commutative.  By step 4 we have an ideal of
$A_0$, $I_1u_1^3$, contained in $Z$, and so $[A_1,I_1u_1^3]=0$, and by lemma 1.3 either $A_1\subseteq
Z(A)_1$ or $I_1u_1^3=0$ for any $u_1\in U_1$. 

If $A_1\subseteq Z_1(A)$ then $A_1^2\subseteq Z$, and
since
$A$ is prime and
$A_1+A_1^2$ is a nonzero ideal contained in $Z(A)$, because $A$ is nontrivial, we deduce that
$A$ is commutative.

If $I_1u_1^3=0$, since $0=I_1u_1^3=(I_1u_1)(u_1^2A)$ and $u_1^2A$ is an ideal of $A$ because
$u_1^2\in Z$ by step 3,  then from the primeness of $A$ either
$I_1u_1=0$ or
$u_1^2=0$ for any
$u_1\in U_1$. But if
$u_1^2=0$ for every $u_1\in U_1$ we get $U_1U_1=0$ because $u_1\circ v_1=0$ for any $u_1,v_1\in U_1$
and if
$I_1u_1=0$ then
$0=I_1(u_1v_1)$ for every $u_1,v_1\in U_1$. From  step 3  and because $A$ is prime we
 obtain that either $U_1U_1=0$ or $I_1=0$. But $I_1=0$ contradicts that $A$ is prime because then
$IA_1=0$ and so $I(A_1+A_1^2)=0$ with $A_1+A_1^2$ a nonzero ideal of $A$. Therefore $U_1U_1=0$ in any
case, when $I_1u_1^3=0$.

\medskip

6. Either $U\subseteq Z$, or $A$ is commutative, or $A$ is a central order in a 4-dimensional simple
superalgebra. We consider
$[v_1,z_1]I$ with
$v_1\in U_1, z_1\in I_1$. It is an ideal of
$A$ by step 1. For any
$u_0\in U_0, v_1
\in U_1$ and $y_1, z_1\in I_1$ we have
$$[u_0,y_1][v_1,z_1]I=[u_0,y_1]v_1z_1I+[u_0,y_1]z_1v_1I$$
{\noindent with
$[u_0,y_1]v_1z_1I=0$ by step 5 and}
$$[u_0,y_1]z_1v_1I=-y_1[u_0,z_1]v_1I+[u_0,y_1z_1]v_1I=0$$
{\noindent by steps 2 and 5. Since $A$ is prime we obtain that either i) $[U_1,I_1]=0$
or ii) $[U_1,I_1]\not= 0$, and then $[U_0,I_1]=0$.}

\medskip

{\parindent = 3em i) If $[U_1,I_1]=0$ then for any $u_1\in U_1$, $u_1I_1$ is a nilpotent ideal of $A_0$ because  by step 5 
$(u_1I_1)(u_1I_1)=u_1^2I_1=0$, and since $A_0$ is semiprime by lemma 1.2, we deduce
that $u_1I_1=0$. But then $u_1I_0A_1\subseteq u_1I_1=0$ and also $u_1I_0A_1^2=0$, that is,
$u_1I_0(A_1+A_1^2)=0$ with $ A_1+A_1^2$ a nonzero ideal of $A$. By the primeness of $A$, $u_1I_0=0$,
and so $u_1I=0$ and $U_1=0$. Therefore $[U,I]=[U_0,I]=[U_0,I_0]$, and $[U_0,I_0]=0$ by step
2, so by lemma 2.3, $U_0\subseteq Z$.}

\medskip

{\parindent= 3em ii) If $[U_1,I_1]\not=0$, then $[U_0,I_1]=0$ and so by step 2 \ $[U_0,I_0]=0$ and from lemma 2.3,
$U_0\subseteq Z$. Also 
$Z\not=0$ and we may localize
$A$ by
$Z$ and consider in
$Z^{-1}A$, the Lie subalgebra $Z^{-1}(ZU)$ and the ideal $Z^{-1}I$, which satisfy the hypothesis of
the theorem. Now we have also that $0\not= Z^{-1}Z$ is a field. By step 1, \  $[U_1,I_1]\subseteq Z$, and
hence}
$$0\not=[Z^{-1}(ZU)_1,Z^{-1}I_1]\subseteq Z^{-1}I_0\cap Z^{-1}Z.$$
{\noindent Therefore $Z^{-1}I$ has invertible elements and so $Z^{-1}I=Z^{-1}A$. But
then  $Z^{-1}(ZU)$ is a Lie ideal of $Z^{-1}A$.  Since $[Z^{-1}(ZU),$
$Z^{-1}(ZU)]=0$ because
$U_0\subseteq Z$ and because of step 5, it follows from theorem 3.2 and its proof in \cite {M} that
either $Z^{-1}(ZU)\subseteq Z^{-1}Z$ or $A$ is a central order in the matrix algebra
$M_{1,1}(Z^{-1}Z)$. In the last case $A$ is a central order in a 4-dimensional simple superalgebra,
and in the first case  $Z^{-1}(ZU)\subseteq Z^{-1}Z$ and we can deduce from the primeness of $A$ that
$U\subseteq Z$.}

\medskip

Therefore in  case a) we have obtained that either $U\subseteq Z$, or $A$ is commutative, or $A$
is a central order in a 4-dimensional simple superalgebra

\bigskip

b) We suppose now that $[T^\prime, I]\not=0$. We recall that $[T^\prime, I]\subseteq T^\prime$.
Consider
$[[T^\prime,I],T^\prime]$.  We claim that $I[[T^\prime,I],T^\prime]I\subseteq T^\prime$. Indeed, let
$x\in T^\prime, y\in [T^\prime,I]$  and  $a\in I$. Since $[T^\prime, I]\subseteq T^\prime$ and
$T^\prime$ is a subring,
$$[x,y]a=[x,ya]-(-1)^{\bar x \bar y}y[x,a]\in T^\prime.$$

{\noindent 
Now, let $b\in I$; we get}
\begin{eqnarray*}
b[x,y]a &=&[b,[x,y]]a + (-1)^{(\bar x+ \bar y) \bar b}[x,y]ba \\ &=& -(-1)^{\bar y ( \bar b +
\ \bar x)}[y, [b,x]]a-(-1)^{\bar b \bar x+ \bar b \bar y}[x,[y,b]]a\\&\ &+(-1)^{(\bar x+ \bar y)\bar
b}[x,y]ba\in T^\prime.
\end{eqnarray*}

{\noindent Therefore, by the primeness of $A$,
$T^\prime$ is dense if $[[T^\prime,I],T^\prime]\not=0$.}

 If $[[T^\prime,I],T^\prime]=0$, then
$$[[[U,I],[U,I]]_0,[[[U,I],[U,I]]_0,I]]=0$$

{\noindent so by theorem 1 in \cite {H3},  $[[[U,I],[U,I ]]_0,I]=0$, and applying now lemma 2.3 we have
$[[U,I],$ $[U,I]]_0 \subseteq Z$. We denote  $V=[[U,I],[U,I]]$ and we have that $V$ satisfies the same conditions as $U$, that is, $V$ is an abelian subgroup of $A$ such that $[V,I]\subseteq V$,  $u\circ v \in Z$ for every $u,v\in V_0$, and $u\circ v=0$ for every $u,v\in V_1$, because $V\subseteq U$. Since $V_0\subseteq Z$ we observe that $V$ has, like $U$ in case a)  steps 1 and 2, the following properties: $[V,I]_0\subseteq V_0\subseteq Z$ and $[V_0,I_0]=0$. From this we can prove steps 3, 4, 5 and 6 in a) exactly in the same way  but now taking $V$ instead of $U$. So we obtain that either $A$ is commutative, or $A$ is a central order in a 4-dimensional simple superalgebra, or $V\subseteq Z$. But if $V\subseteq Z$ we can apply case a) and we obtain  that either $U\subseteq Z$, or
$A$ is commutative, or $A$ is a central order in a 4-dimensional simple superalgebra.}

It remains to consider the case when $T^\prime$ is dense in $A$.  We denote by
$J=I[[T^\prime,I],T^\prime]I$ and so $J\subseteq T^\prime$. From the definition of
$T$  and because $T^\prime \subseteq T$ we know that $[T^\prime, A]\subseteq [U,I]$, and therefore
$[J,A]\subseteq [U,I]\subseteq U$. By hypothesis $u\circ v \in Z$ for
any $u,v\in U_0$, so $u\circ v\in Z$ for any $u, v\in [J,A]_0$.

We assume first that $u\circ v=0$ for any $u,v\in  [J,A]_0$. Then $1/2 (u\circ u)=u^2=0$ for any
$u\in [J,A]_0$ and since $A_0$ is semiprime by lemma 1.2, we can apply lemma 1 in \cite {La-M} and we
have $[J,A]_0=0$. Therefore $[J,A]=[J,A]_1$ and then $[J,A]$ is a Lie ideal of $A$ such that
$[[J,A],[J,A]]=0$. From theorem 3.2 and its proof in \cite {M} it follows that either $[J,A]\subseteq
Z$ or $A$ is a central order in a 4-dimensional matrix superalgebra. If $[J,A]\subseteq Z$, since
$[J,A]=[J,A]_1$, we get $[J,A]=0$ and now by lemma 2.3,  $J\subseteq Z$, and so $A$ is commutative.

Suppose now that there exist $u,v\in [J,A]_0$ such that $u\circ v\not= 0$. Then $Z\not= 0$, and we
may form the localization $Z^{-1}A$. Since $[J,A]\subseteq [U,I]\subseteq U$ we have
$[Z^{-1}J,Z^{-1}A]\subseteq [Z^{-1}(ZU), Z^{-1}I]\subseteq Z^{-1} (ZU)$, and so from the hypothesis
of the theorem for any $u,v\in [Z^{-1}J, Z^{-1}A]_0$ we get $u\circ v\in Z^{-1} Z \cap Z^{-1}J$. But
$Z^{-1}Z$ is a field and so $Z^{-1}J$ has some invertible element forcing $Z^{-1}J=Z^{-1}A$. Therefore
$[Z^{-1}J, Z^{-1}A] = [Z^{-1}A, Z^{-1}A]\subseteq Z^{-1}(ZU)$ and again by the hypothesis of
the theorem it follows that $[Z^{-1} A, Z^{-1} A]_1 \circ [Z^{-1} A, Z^{-1} A]_1=0$. We apply now
lemma 2.6 in \cite {M} and we obtain that $Z^{-1}A$ is commutative (superalgebras of the type (b) and
(c) in the lemma do not satisfy the condition $[Z^{-1} A, Z^{-1} A]_1 \circ [Z^{-1} A, Z^{-1}
A]_1=0$), and so $A$ is commutative. This finishes the proof. 
\end{proof}

Next we consider the cases when $P^*= P$ and the involution on $A/P$ is of the second kind or of the first kind.

\begin{lem}
Let $A$ be a prime superalgebra with a superinvolution $*$ of the second kind. Let $U$ be a Lie ideal
of $K$ such that $u\circ v \in Z$ for every $u,v \in U_0$, and $u\circ v=0$ for every $u,v\in U_1$.
Then either $U\subseteq Z$ or $A$ satisfies $S(2)$. 
\end{lem}

\begin{proof}[Proof:]
If $*$ is of the second kind we know that $Z_H=\{x\in Z : x^*=x\}\not= Z$. We may localize $A$ by
$V$ and replace $U$ by $V^{-1}(Z_HU)$ and $A$ by $V^{-1}A$. The hypothesis remains unchanged, so we
keep for this superalgebra the same notation $A$, and now
$Z$ is a field. Let
$0\not= t\in Z_K$. Then $H=tK$ and $A=tK+K$. It follows that
$[ZU,A]\subseteq ZU$,  $u\circ v \in Z$ for every $u,v \in ZU_0$, and $u\circ v=0$ for
every $u,v\in ZU_1$. By theorem 2.4, either $ZU\subseteq Z$, which implies that $U\subseteq Z$, or $A$
satisfies $S(2)$.
\end{proof}

\begin{lem}
Let $A$ be a prime superalgebra with a superinvolution $*$ of the first kind. Let $U$ be a Lie ideal
of $K$ such that $u\circ v \in Z$ for every $u,v \in U_0$, and $u\circ v=0$ for every $u,v\in U_1$.
Then either $U\subseteq Z$ or $A$ satisfies $S(4)$. 
\end{lem}

\begin{proof}[Proof:]
If $u^2=0$ for every $u\in U_0$, applying theorem 3.3 in \cite {Go-L-S} we obtain that $U=0$. Suppose
then that $u^2\not= 0$ for some $u\in U_0$. By theorem 3.4 in \cite {Go-L-S} we get that either
$U\subseteq Z$ or $A$ is a central order in a Clifford algebra with either 2 or 4 generators.

\end{proof}

\bigskip

Combining the above results we obtain

\bigskip

\begin{theo}
Let $A$ be a semiprime superalgebra and $U$ a Lie ideal of $K$ with $u\circ v \in Z$ for every $u,v
\in U_0$, and $u\circ v=0$ for every $u,v\in U_1$. Then $A$ is the subdirect sum of two semiprime
homomorphic images $A^{\prime}$, $A^{\prime \prime}$, such that $A^\prime$ satisfies $S(4)$ and the
image of $U$ in $A^{\prime \prime}$ is central.
\end{theo}

\begin{proof}[Proof:]
Let $T^\prime=\{ P: P$ is a prime ideal of $A$ such that $A/P$ satisfies $S(4)\}$ and let
$T^{\prime \prime}= \{ P : P$ is a prime ideal of $A$ such that the image of $U$ in $A/P$ is
central$\}$. 

If we consider $P$ a prime ideal of $A$ such that $P^*\not= P$ we know from theorem 2.8
that either $A/P$  is a central order in a simple superalgebra at most 4-dimensional
over its center, or $(U+P)/P$ is central. If we consider $P$ a prime ideal of
$A$ such that $P^*=P$, it follows from lemmata 2.5, 2.6 that either $A/P$  is a central
order in a simple superalgebra at most 16-dimensional over its center, or the image of $U$ in $A/P$ is
central. 

So every prime ideal of $A$ belongs either $T^\prime$ or $T^{\prime \prime}$. Then
$A^\prime$ is obtained by taking the quotient of $A$ by the intersection of all the prime ideals in
$T^\prime$, and $A^{\prime \prime}$ is obtained by taking the quotient of $A$ by the intersection of
all the prime ideals in $T^{\prime \prime}$. This proves the theorem.
\end{proof}

\bigskip

We finally arrive at the main theorem on the Lie structure of $K$.

\bigskip

\begin{theo}
Let $A$ be a semiprime superalgebra with superinvolution $*$, and let $U$ be a Lie ideal of $K$. Then
either $A$ is a subdirect sum of two semiprime homomorphic images $A^\prime$, $A^{\prime \prime}$,
with $A^\prime$ satisfying $S(4)$ and the image of $U$ in $A^{\prime \prime}$ being central, or $U\supseteq
[J\cap K, K]\not= 0$ for some ideal $J$ of $A$.
\end{theo}

\begin{proof}[Proof:]
From lemmata 2.1 and 2.2 we know that either $U$ is dense in $A$, and so there exist a
nonzero ideal $J$ such that $J\subseteq \bar U$, or $u\circ v \in Z$ for every $u,v \in U_0$, and
$u\circ v=0$ for every $u,v\in U_1$. In the second case we obtain by theorem 2.7 the first part of
the theorem. So suppose that $J\subseteq \bar U$.  

The identity
$$[xy,z]=[x,yz]+(-1)^{\bar x \bar y+\bar x \bar z}[y,zx]$$

{\noindent can be used to show that $[\bar U,A]=[U,A]$. Hence $[J\cap K, K]\subseteq [\bar U,
A]=[U,A]=[U,H]+[U,K]$. But $[U,H]\subseteq H$, and $[U,K]\subseteq K$, so $[J\cap K, K]\subseteq [U,K]
\subseteq U$.}

Finally, suppose that $[J\cap K, K]=0$, then $[u\circ v,w]=0$ for any $u,v,w \in J\cap K$ because
$[uv,w]=u[v,w]+(-1)^{\bar v
\bar w}[u,w]v=0$. So by lemmata 2.1, 2.2 and theorem 2.7 it follows
that for each prime image, $A/P$, of $A$ either  its center contains
$((J\cap K)+P)/P$, or $A/P$ is a central order in a simple superalgebra at most 16-dimensional over
its center.

We claim that if the image of $J \cap K$ in $A/P$ for some prime ideal $P$ of $A$ is central, then $A$ is as described in the first part of the conclusion of the theorem. 

Let $P$ be a prime ideal such that $P^*\not= P$. If $(J+P)/P\not=0$,
then since $A/P$ is a prime superalgebra we get $((J\cap P^*) +P)/P \not=0$, and so we have $((J\cap
P^*) +P)/P\subseteq ((J\cap K)+P)/P \subseteq Z_0(A/P)$, that is, $A/P$ is commutative. So $A/P$ is
commutative unless $J\subseteq P$. And if $J\subseteq P$, then by the proof of lemma 2.1  we know
that  $A[u\circ v, w]A[u\circ v,w]A\subseteq P$ for any $u, v, w \in U$, and because  $P$
is a prime ideal we deduce that $[u\circ v,w]\in P$ for any $u, v, w \in U$. But now by Lemma 2.2 and
since
$[u\circ v,w]+P=0$ for any $u,v,w \in U$,  it follows that
$A/P$ satisfies the conditions $u\circ v\in Z$ for any $u,v\in ((U+P)/P)_0$ and
$u\circ v=0$ for any $u,v \in ((U+P)/P)_1$. By theorem 2.4 we obtain that either
$(U+P)/P\subseteq Z_0(A/P)$, or $A/P$ satisfies $S(4)$.

 And if $P$ is a prime ideal such that $P^*=P$ then $A/P$ has a superinvolution
induced by * and
$K(A/P)=(K+P)/P$.  In this case if $((J\cap K) +P)/P=0$ we get $(J+P)/P\subseteq (H+P)/P=H(A/P)$, and
therefore $(J+P)/P$ is supercommutative. But then for any $a,b \in A/P$ and $y,z \in (J+P)/P$ it
follows that

\begin{eqnarray*}
yabz&=&(-1)^{(\bar b + \bar z)(\bar y + \bar a)}(bz)(ya)=(-1)^{\bar b (\bar y +
\bar a)} b(ya)z\\
&=& (-1)^{\bar b \bar y + \bar b \bar a+ (\bar a +\bar z)\bar
y} b(az)y= (-1)^{\bar b \bar a} ybaz,
\end{eqnarray*}

{\noindent and since $A/P$ is prime $ab = (-1)^{\bar a \bar b} ba$, that is, $A/P$ is
supercommutative. Now from lemma 1.9 in \cite {M},  $A/P$ is a central order in a simple superalgebra
at most $4$-dimensional over its center. And if
$((J\cap K)+P)/P\not=0$ then
$Z_0(A/P)\not=0$,  so by localizing at
$V=(Z_0(A/P)\cap H(A/P))-\{0\}$  we can suppose that $Z_0(A/P)$ is a field, which we denote by
$Z$. We will replace $V^{-1}(A/P)$ by $A/P$ and $V^{-1}((J+P)/P)$ by $(J+P)/P$. Then if
$0\not= t\in ((J\cap K)+P)/P$ we have  $tH=K$ with
$H=H(A/P), K=K(A/P)$,  so $K=tH\subseteq K\cap J\subseteq Z$, and also $tH=K\subseteq Z$ and
$H\subseteq t^{-1}Z\subseteq Z$. Therefore $A/P$ is a field.}

\end{proof}

Finally we have

\bigskip

\begin{coro}
Let $A$ be a semiprime superalgebra with superinvolution $*$, and let $U$ be a Lie ideal of $K$. Then
either $[J\cap K,K]\subseteq U$ where $J$ is a nonzero ideal of $A$ or there exists a semiprime ideal
$T$ of $A$ such that $A/Ann T$ satisfies $S(4)$ and $(U+T)/T\subseteq Z_0(A/T)$.
\end{coro}

\begin{proof}[Proof:]
By theorem 2.8 we have that either the first conclusion holds, or, for each prime ideal $P$ of $A$,
either $A/P$ satisfies $S(4)$ or $(U+P)/P \subseteq Z_0(A/P)$. Let $T$ be the intersection of the
prime ideals $P$ of $A$ such that $(U+P)/P\subseteq Z_0(A/P)$. Then
$Ann T$ contains the intersection of those prime ideals $P$ such that $A/P$ satisfies $S(4)$. So we
get that
$A/Ann T$ satisfies $S(4)$, and this proves the result.
\end{proof}

\end{document}